\documentclass[11pt,reqno]{amsart}
\usepackage{amsmath,amsfonts,amssymb,amsthm,amscd}

\makeindex

\usepackage{amssymb}

\usepackage[latin1]{inputenc}

\usepackage{epsfig,epsf}
\usepackage{enumerate}
\usepackage{indentfirst}
\usepackage{euscript}
\usepackage{graphicx}
\usepackage{amsmath}
\usepackage{amsfonts}
\usepackage{amssymb}
\makeindex

 \theoremstyle{plain}

\newtheorem{proposition}{Proposition}

\newtheorem{corollary}{Corollary}
\newtheorem{remark}{Remark}
\theoremstyle{definition}

\usepackage{amsmath,amsfonts,  amssymb,amsthm,amscd}
\makeindex

\usepackage{amssymb}
\usepackage{graphics}

\usepackage[latin1]{inputenc}
\usepackage{epsf}

\usepackage{ epsfig,epsf}
\usepackage{enumerate}
\usepackage{indentfirst}
\usepackage{euscript}
\usepackage{graphicx}
\usepackage{amsmath}
\usepackage{amsfonts}
\usepackage{amssymb}
\makeindex

 \theoremstyle{plain}

\title[Curvature Lines around  Umbilic Curves]{On the Patterns of
Principal  Curvature Lines around  a Curve of Umbilic Points}

 \author[R. Garcia]{Ronaldo Garcia$^1$}

 \author[J. Sotomayor]{Jorge Sotomayor$^2$}

\thanks{ The  authors are fellows of  CNPq.
 This work was done under the project  CNPq/PADCT 620029/2004-8. The second
  author acknowledges the help of P. Tonelli and E. Harle in
   translating  (Carath\'eodory 1935). \\ Correspondence to: ragarcia@mat.ufg.br / sotp@ime.usp.br}

\begin{document}
 \maketitle

 {\small \address{ \centerline{  $^1$ Instituto de Matem\'{a}tica e
Estat\'{\i}stica  }
 \centerline{Universidade Federal de Goi\'as}
  \centerline{ Caixa Postal 131}
\centerline{74001-970  Goi\^ania, GO, Brasil}
 \centerline{\\}
\centerline{\\} \centerline{$^2$ Instituto de Matem\'{a}tica e
Estat\'{\i}stica}
  \centerline{Universidade de S\~{a}o Paulo}
   \centerline{Rua do Mat\~{a}o 1010,
Cidade Universit\'{a}ria}
   \centerline{05508-090  S\~{a}o Paulo, SP, Brasil}
  }}

 \begin{abstract}
 In this paper is studied the behavior  of principal curvature lines near a curve
of umbilic points of a smooth
   surface.
   \end{abstract}


\section{Introduction}

The study of umbilic points on surfaces and the patterns of
principal curvature
 lines around them has attracted the attention of generation of
mathematicians
  among whom can be named Monge, Darboux and  Carath\'eodory.
  One aspect --concerning  isolated umbilics-- of the
   contributions of these authors has been elaborated and extended
   in several direction by Garcia,  Sotomayor
   and Gutierrez,  among others. See
   (Gutierrez and Sotomayor, 1982, 1991, 1998), (Garcia and Sotomayor,
   1997, 2000) and  (Garcia et al. 2000, 2004)
    where additional references can be found.

  In  (Carath\'eodory 1935) Carath\'eodory  mentioned  the interest of
  non isolated umbilics in  generic surfaces pertinent to
   Geometric Optics.   In a remarkably concise study he    established
that any
local analytic  regular arc of curve in  $\mathbb R^3$ is a curve
of umbilic points of  a piece of  analytic surface. In some cases
he also determined the patterns  of behavior of principal
curvature lines near the curve of umbilic points.

In the present
 paper will be  performed an analytic,  explicit
 and
 constructive study
 of umbilic curves and of the simplest patterns for their
 neighboring principal curvature lines, that holds  also  for   smooth
  curves and surfaces.  A comparison of results of this work  with
  those of
  Carath\'eodory (Carath\'eodory 1935)  is attempted in section \ref{sec:cr}.

\section{Preliminaries} \label{sec:pr}

Let $c:[0,l]\to \mathbb R^3$ be a regular curve parametrized by
arc length $u$ contained in a regular smooth surface $\mathbb M$,
which is
 oriented by the once for all given  positive unit normal vector field
$N$.

Let $T\circ  c =c^\prime$. According  to (Spivak 1980),  the
Darboux frame $\{T, N\wedge T,N\}$ along  $c$ satisfies the
following system of differential equations:

\begin{equation}\label{eq:da}
\aligned T^\prime &= k_g N\wedge T +k_n N\\
(N\wedge T)^\prime &=-k_g T +\tau_g N\\
N^\prime  &=-k_n T - \tau_g (N\wedge T)\endaligned
\end{equation}

\noindent where $k_n$ is the {\it normal curvature}, $k_g$ is the
{\it geodesic curvature} and $\tau_g$ is the {\it geodesic
torsion} of the  curve $c$.

\begin{proposition} \label{prop:ca} Let $c: [0,l]\rightarrow
\mathbb M$ be a regular arc length parametrization of a curve of
umbilic points,  such that $\{ T,N\wedge T, N\}$ is a positive
frame of ${\mathbb R}^3$. Then the expression
\begin{equation}\label{eq:ca} \alpha (u,v)= c(u) + v   (N\wedge T)(u)
+ [\frac 12
  k (u)v^2 +\frac 16 a(u) v^3+ \frac 1{24} b(u)v^4+ h.o.t]N(u),
\end{equation}

\noindent where  $k(u)=k_n(c(u),T)=k_n(c(u),N\wedge T) $ is the
normal curvature of $\mathbb M$ in the directions $T$ and $N\wedge
T$, defines a local $C^\infty $ chart
 in a small tubular
neighborhood of $c$. Moreover $\tau_g(u)=0$.
\end{proposition}

\begin{proof} This parametrization, in the case where $c$ is   a
principal curvature line, was first introduced by  Gutierrez and
Sotomayor in (Gutierrez and Sotomayor, 1982).
 In the
present case, $c$ being a regular curve of umbilic points,  by the
Implicit Function Theorem it follows  that the principal
curvatures $k_1\leq k_2$ restricted to $c$ are also $C^\infty$.
Also, it follows that $k_n(c(u), T(u))=k_n(c(u), (N\wedge T)(u))$
and in consequence $\tau_g(c(u),v)=0$ for any $v\in
T_{c(u)}\mathbb M$.
\end{proof}

\section{ First and Second Fundamental Forms}

In the chart $\alpha$ in the equation (\ref{eq:ca}) the positive
unit normal vector vector field  is given by
$$ N=
\frac{\alpha_u\wedge \alpha_v}{|\alpha_u\wedge \alpha_v|}.$$
So it follows that,
$$\aligned   N(u,v)=&-[ \frac 12k^\prime v^2+\frac
16(a^\prime+3k^\prime
k_g)v^3+0(v^4)] T(u)\\ -&[k v+\frac 12 a(u) v^2+ \frac 16( b(u) -3
k
^3) v^3+0(v^4)](N\wedge T)(u)\\
+&[1-\frac 12k^2 v^2+\frac 12k^\prime k v^3+O(v^4)]N(u)\endaligned
$$

 The coefficients of the {\it first} and
{\it second fundamental forms} in the chart $\alpha$ are given by:
$E=<\alpha_u,\alpha_u>, \;\;  F=<\alpha_u,\alpha_v>,$\- \;
$G=<\alpha_v,\alpha_v>,$ \- \;\;
$e=<\alpha_{uu},N>,\;\;f=<\alpha_{uv},N>$ and $g=<\alpha_{vv},N>.$
Therefore,

\begin{equation}
\aligned
E(u,v)=& 1-2k_g v+(k_g ^2-k ^2)v^2+\frac 16(6k_gk^2-2ka(u))v^3+
O(v^4)\\
F(u,v)=&\frac 12k^\prime kv^3+O(v^4)\\
G(u,v)=&  1+k ^2v^2+ka(u)v^3+O(v^4) \\
e(u,v)=& k -2k_g k v+\frac 12
 (2 k
 k_g^2-k_g a(u)-2k^3+k^{\prime\prime} )v^2\\
 +&   \frac 16[ \;a^{\prime\prime}+k_g(9k^3-b(u))+ (3k_g^2-k^2)a(u)  \\
 +&  3k^\prime( k_g^\prime + k^2)\; ]v^3
                    +O(v^4)\\
 f(u,v)=& k^\prime  v+\frac 12(k_g k^\prime
 +a^\prime )v^2+\frac 16 (k_g a^\prime+3k^\prime k_g^2+b^\prime)v^3+
O(v^4)\\
 g(u,v)=& k +a(u)v+\frac 12(b(u)-k ^3)v^2-\frac 12
k^2(a(u)-k^\prime)v^3+O(v^4)
 \endaligned
 \end{equation}

 The {\it Mean} and {\it Gauss}  curvatures in the chart
$\alpha$  are given by:

 \begin{equation}\label{eq:hk}
 \aligned
{\mathcal H}=&k +\frac 12a(u)v+\frac 14( b(u)+k^{\prime\prime}-3
k^3- k_g a(u))v^2+O(v^3)\\
{\mathcal K}=&k^2+k a(u)v+\frac 12(-k_g k a(u)-3 k^4+k
k^{\prime\prime}+k b(u)-2k^{\prime\prime})v^2+O(v^3)
\endaligned
 \end{equation}

 According to
 (Spivak 1980) the differential equation of curvature lines
in the chart $\alpha $  is given by

 \begin{equation}\label{eq:cl}
  \aligned (F&g-Gf)dv^2+(Eg-Ge)dudv+(Ef-Fe)du^2=\\
  =&Ldv^2+Mdvdu+Ndu^2=0;  \;\;\;\;  \;{\text  where} \\
   L=&-[ k^\prime v+\frac 12(k_g k^\prime +a^\prime)
 v^2+ \frac 16(k_g a^\prime+ 3k^\prime k_g^2+ b^\prime+
3k^2k^\prime)v^3+ O(v^4) ]  \\
 M=&\;\;  a(u)v+\frac 12[b(u)-3k^3-k^{\prime\prime}-3k_g a(u)]v^2\\
+&  \frac 16[ 15k^3k_g - 3k_g^\prime k^\prime +
(3k_g^2-16k^2)a(u)-a^{\prime\prime} -5k_g b(u) ]v^3+O(v^4) \\
 N=& \;k^\prime v+\frac 12(a^\prime-3k_g k^\prime )
 v^2+ \frac 16( 3k^\prime k_g^2 -9k^2k^\prime -5k_g
a^\prime+b^\prime)v^3+
 O(v^4)
 \endaligned\end{equation}

\section{Principal Configurations near an Umbilic Curve}

\begin{figure}[htbp]
\begin{center}
\includegraphics[angle=0, width=10.5cm]{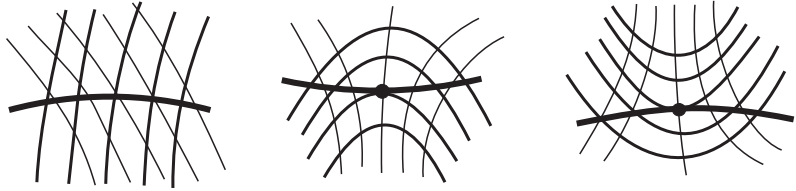}
\end{center}
\caption{ \label{fig:1} Principal curvature lines near an umbilic
curve: transversal case, left,  and tangential  case, center and
right.}
\end{figure}

\begin{proposition}\label{prop:1} Suppose that $\nabla {\mathcal
H}(u,0)=(k^\prime, a(u)/2)$ is not zero at a point $u_0$.
 Then the
principal foliations near the point  $c(u_0))$ of the curve are as
follows.

\begin{itemize}
\item[i)] If $k^\prime(u_0)\ne 0$ then both principal foliations
are transversal to the curve of umbilic points. See Fig.
\ref{fig:1}, left.

\item[ii)] If $k^\prime(u_0)= 0$,  $k^{\prime\prime}(u_0)\ne  0$
and $a(u_0)\ne 0$, then one principal foliation is transversal to
$c$ and the other foliation has quadratic contact with the curve
$c$ at the point $c(u_0)$. See Fig. \ref{fig:1}, center and right.
\end{itemize}
\end{proposition}

\begin{proof} After division by $v$,  the implicit differential equation
of curvature lines is given by:

 \begin{equation}\label{eq:cl1}
 \aligned   {\mathcal F}(u,v,[du:dv])= -&[ k^\prime  +\frac 12(k_g
k^\prime +a^\prime)
 v +  O(v^2) ] dv^2\\
+&[ a(u) +\frac 12(b(u)-3k^3-k^{\prime\prime}-3k_g a(u))v
 +O(v^2)]dudv\\
+&[k^\prime  +\frac 12(a^\prime-3k_g k^\prime )
 v + O(v^2)]du^2=0, \;\;\;
 \endaligned\end{equation}

Consider the Lie - Cartan line field $X=({\mathcal F}_p,
p{\mathcal F}_p, -({\mathcal F}_u+p{\mathcal F}_v))$ where
${\mathcal F}=Lp^2+Mp+N=0, \; p=dv/du$. It follows that ${\mathcal
F}(u,0,p)=k^\prime(1-p^2)+a(u)p$ and so  when $   k^\prime\ne 0$,
${\mathcal F}(u,v,p)=0$  is a regular surface near $v=0$ and has
two connected components containing the points $(u,0,p_{\pm})$,
where $p_{\pm}=[a(u)\pm \sqrt{a^2+4(k^\prime)^2}
\;]/(2k^\prime)\ne 0$. Clearly the Lie-Cartan line field field $X$
is regular at $(0,0,p_{\pm})$ and so the configuration of
principal lines is as shown in Fig. \ref{fig:1}, left.

 Suppose for simplicity that at $u=0$ it holds that
$k^\prime(0)=0$. By hypothesis, we have that $a(0)=a_0\ne 0$. Then,
after
division by $v$,  it follows that the implicit differential
equation of principal curvature lines is given by:

$$\aligned {\mathcal F} =&-[k^{\prime\prime}(0)u+\frac 12 a^\prime(0)
v+
\cdots]dv^2\\
+&[a_0+ a^\prime(0)u+ (b(0)-3k(0)^3-k^{\prime\prime}-3k_g(0)a_0)
v+\cdots]dvdu\\
+ & [k^{\prime\prime}(0)u+\frac 12 a^\prime(0) v+
\cdots]du^2=0\endaligned$$

 Therefore the  directions defined by ${\mathcal
 F}(u,v,[du:dv])=0$ are $p=0$ and $q=0$ where $p=dv/du$ and
 $q=du/dv$. Therefore one  direction is tangent to $ c$
 and the other is orthogonal to $c$.

By the conditions imposed,  one foliation is orthogonal to, and
the other has quadratic contact with,  the umbilic curve  at
$c(0)$.

In fact, ${\mathcal F}_p(0,0,0)=a_0\ne 0$ and it follows that the
integral curve of the Lie - Cartan line field $X=({\mathcal F}_p,
p{\mathcal F}_p, -({\mathcal F}_u+p{\mathcal F}_v))$ passing
through $(0,0,0)$ is given by:
$$\aligned u(t)=& \; a_0t+\cdots\\
v(t)=&-\frac 12 k^{\prime\prime}(0) a_0 t^2 +\cdots \\
p(t)=&-k^{\prime\prime}(0)t+\dots
\endaligned$$

The  curve $(u(t),v(t))$  has  quadratic contact with the line
$v=0$, provided  $k^{\prime\prime}(0) a_0\ne 0$. Eliminating $t$
it follows that $v=-\frac{k^{\prime\prime}}{2a_0} u^2+\cdots$.

 To describe the
other principal foliation near $(0,0,0)$ it is convenient to
consider the Lie - Cartan line field $Y=(q{\mathcal F}_q,
{\mathcal F}_q, - ( q {\mathcal F}_u  + {\mathcal F}_v))$ with
$q=du/dv$. It is clear that this  principal foliation is
transversal to $c$ in a neighborhood of $c(0)$. Therefore the
principal configuration is as shown in Fig \ref{fig:1}, center and
right. This ends the proof.
\end{proof}

\begin{proposition}\label{prop:11} Suppose that
$k^\prime(0)=a( 0)=0 $, $a^\prime(0)k^{\prime\prime}(0)\ne 0$,  at
the point   $c(0)$ of a regular curve $c$ of umbilic points. Let
$A:= -2k^{\prime\prime}(0)/a^\prime ( 0)\ne 0$ and $B:=
[b(0)-3k(0)^3-k^{\prime\prime}(0)]/a^\prime (0)$. Let $\Delta$ and
$\delta$ be defined by
$$\Delta= -4 A^4+12 B A^3-(36+12 B^2) A^2+( 4 B^3+72 B) A-9
B^2-108; \; \; \delta =2- AB.$$

Then the principal foliations at this point are as follows.

\begin{itemize}

\item[i)] If   $ \delta <0 $ and $\Delta <0$ then $0$ is
topologically equivalent to  a Darbouxian umbilic of type $D_1$,
through  which the umbilic curve is adjoined  transversally  to
the
 separatrices.
 See Fig.
\ref{fig:dl} left.

\item[ii)] If   $ \delta  <0 $ and $\Delta>0$ then $0$ is
topologically equivalent to a Darbouxian umbilic of type $D_2$,
through which the  umbilic
 curve is adjoined,  on the interior of the parabolic sectors,
transversally  to the
 separatrices and to  the nodal central line.
See     Fig. \ref{fig:dl} center.

\item[iii)] If $  \delta >0$ then  $0$ is topologically  a
Darbouxian umbilic of type $D_3$,   through which the umbilic
curve is adjoined transversally  to the
 separatrices.
See   Fig. \ref{fig:dl} right.
\end{itemize}

\begin{figure}[htbp]
\begin{center}
\includegraphics[angle=0, width=10.5cm]{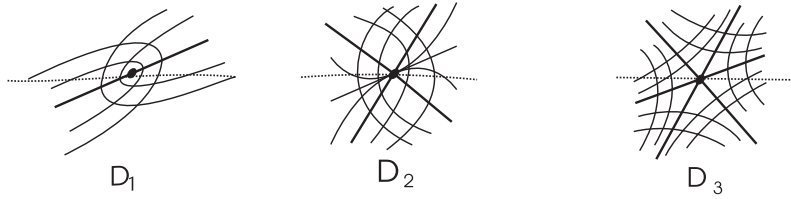}
\end{center}
\caption{ \label{fig:dl} Principal curvature lines near a
Darbouxian-like point on an  umbilic curve, dotted.}
\end{figure}
\end{proposition}

\begin{proof} After division by $v$, the first order terms of
differential equation of curvature lines, see equation
(\ref{eq:cl1}), is given by:
$$\aligned {\mathcal F}(u,v,[du:dv])=
&-(Au+v+O_1(2))dv^2+(2u+Bv+O_2(2))dudv\\
+&(Au+v+O_3(2))du^2. \endaligned$$

 Consider the Lie-Cartan line field
$X=({\mathcal F}_p, p{\mathcal F}_p, -({\mathcal F}_u+p{\mathcal
F}_v))$, where $p=dv/du$. The singular points of  $X$  along the
projective line -- represented by the $p$ axis  \; -- are defined
by the cubic  equation $R(p)=p^3+(A-B) p^2-3 p-A=0$. The condition
$\Delta \ne 0$ means that $R$ has no multiple roots and
consequently  the singularities of $X$ restricted to the
projective line are hyperbolic.    The linear part of $X$ along
the axis $p$ is given by:

$$\aligned DX(0,0,p)= & \left( \begin{matrix}    2-2 A p & B-2 p  & 0
&\\
(2-2 A p ) p &(B-2 p ) p&0 \\  0&0& 3p^2+2( A-B)p-3
\end{matrix}\right)\endaligned$$

The non vanishing eigenvalues of $DX(0,0,p)$ are
$$ \lambda_1=
2+(B-2A)p-2p^2, \;\;\; \lambda_2=3 p^2+2 (A-B) p -3.$$ The
resultant between $R$ and $\lambda_2=R^\prime$ is exactly
$\Delta\ne 0$, while the resultant between $R$ and $\lambda_1$ is
$(2-AB)(16+(2a-b)^2)=\delta (16+(2a-b)^2)\ne 0$.

Therefore in the hypothesis considered all the singular points of
$X$ are hyperbolic. In the chart $(u,v,q=du/dv)$ the Lie-Cartan
line field $Z=(q{\mathcal F}_q, {\mathcal F}_q,-( q{\mathcal F}_u+
{\mathcal F}_v))$ is regular at $0$.

Now consider the vector field $Y=(  Au+v+ O_1(2), -u-\frac B2
v-O_2(2))$. By properties of the index
 of umbilic points,
see (Gutierrez and Sotomayor, 1998), (Hopf 1979) and (Spivak
1980),  it follows that $Ind(Y,0)=-2 Ind({\mathcal F}_i,0)$, where
${\mathcal F}_i$, $(i=1,2)$ denotes the principal foliations. As
$det(DY(0))= (2-AB)/2=\delta/2$,  it follows that $Ind(Y,0)=1$
provided $\delta
>0$,  and that $Ind(Y,0) = -1$ when $\delta <0$. Therefore, $
Ind({\mathcal F}_i,0)=\pm 1/2$.

In the hypothesis of item i), $\Delta<0$  and $\delta <0$,   the
field $X$ has only one hyperbolic saddle point and the index of
the principal foliations is $1/2$ and so they define topologically
a Darbouxian umbilic $D_1$, with one separatrix approaching  the
umbilic point and one hyperbolic sector for each principal
foliation.

In the hypothesis of item ii), $\Delta>0$  and $\delta <0$,   the
field $X$ has two   hyperbolic saddle points and one hyperbolic
node. In this case the index of the principal foliations is $1/2$
and so is topologically a Darbouxian umbilic $D_2$, with  two
separatrices approaching  the umbilic point,  one hyperbolic
sector and one parabolic sector for each principal foliation.

In the hypothesis of item iii), $\Delta>0$  and $\delta >0$,   the
field $X$ has three   hyperbolic saddle points  and the index of
the principal foliations is $-1/2$ and they define  topologically
a Darbouxian umbilic $D_3$, with  three  separatrices approaching
to the umbilic point and three  hyperbolic sectors  for each
principal foliation.

Figure \ref{fig:darbp} shows the behavior of $X$ near the
projective line and the blowing down of the integral curves.

Clearly when $A\ne 0$ the umbilic separatrices are transversal to
the curve of umbilic points which is  parametrized by $v=0$.

\begin{figure}[htbp]
\begin{center}
\includegraphics[angle=0, width=8.5cm]{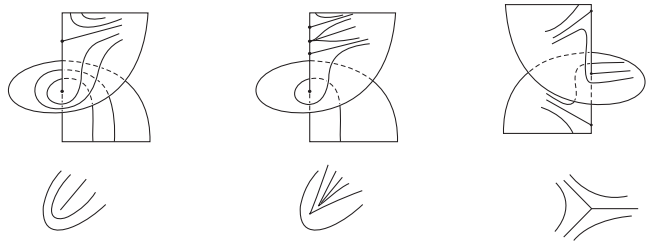}
\end{center}
\caption{ \label{fig:darbp} Behavior of $X$   in a neighborhood of
the projective line and projection of the integral curves }
\end{figure}
\end{proof}

\section{Spherical and Planar Umbilic Curves }

In the previous  sections have been studied a sample of  the most
generic situations, under the restriction on  the surface  of
having an umbilic curve. Below will be considered the case where
$k$ is a constant which implies the additional constrain that the
umbilic curve be spherical or planar, a case also partially
considered in (Carath\'eodory 1935). Under this double imposition
the simplest patterns of principal curvature lines are analyzed in
what follows.

\begin{figure}[htbp]
\begin{center}
\includegraphics[angle=0, width=10.5cm]{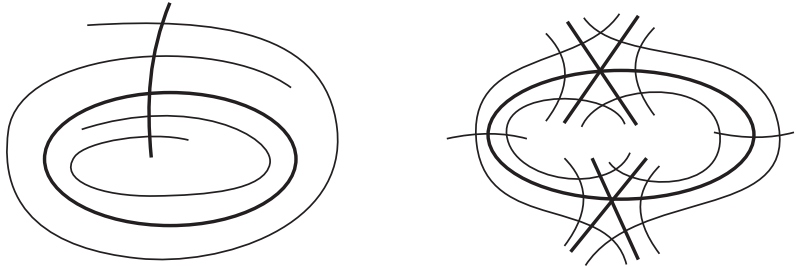}
\end{center}
\caption{ \label{fig:2} Curvature lines near a spherical  umbilic
curve  }
\end{figure}

\begin{proposition}\label{prop:4} Let $c$ be a regular closed spherical
 or planar
curve. Suppose that $c$ is a regular curve of umbilic points on a
smooth surface.  Then  the principal foliations near the curve are
as follows.

\begin{itemize}
\item[i)] If ${\mathcal H}_v(u,0)=a(u)/2 \ne 0$ and  $a(u)>0$ for
definiteness,  then one principal foliation is transversal to the
curve $c$ of umbilic points.

The other foliation defines a first return map (holonomy) $\pi$
along the oriented umbilic curve $c$, with first derivative
$\pi^{\prime}=1$ and second derivative given by a positive
multiple of
 $$\int_0^l k_g(u)
\frac{a^\prime(u)}{a(u)^{\frac 32}}du.$$
  When the above integral is non zero the principal lines  spiral
towards or away from $c$,
   depending on their side relative to
$c$.

\item[ii)] If $a(u)$ has only transversal zeros, near them    the
principal foliations have the topological behavior of a Darbouxian
umbilic point $D_3$ at
 which a separatrix has been replaced  with the umbilic curve.
See Fig. \ref{fig:2}
\end{itemize}
\end{proposition}

\begin{proof} When $c$ is a spherical or planar curve it follows that
$k(u)=k$. In the planar case $k=0$ and $k_g$ is the curvature of
plane curves. In this case the differential equation of principal
curvature lines, see equation (\ref{eq:cl1}), after division by
$v$ is given by:
 \begin{equation}\label{eq:cle}
 \aligned    -&[   \frac 12 a^\prime(u)
 v + \frac 16( k_g a^\prime  + b^\prime )v^2+ O(v^3) ] dv^2\\
+&[ a(u) +\frac 12(b(u)-3k^3 -3k_g a(u))v \\
+&  \frac 16( 15k^3k_g   + (3k_g^2-16k^2)a(u)-a^{\prime\prime}
 -5k_g b(u) )v^2+O(v^4)]dudv\\
+&[ \frac 12 a^\prime(u)
 v + \frac 16(   b^\prime-5k_g a^\prime)v^2+ O(v^3)]du^2=0
 \endaligned\end{equation}

The curve $v=0$ is a solution of equation (\ref{eq:cle}). As
$a(u)\ne 0$ by hypothesis it follows that $v=0$ is a periodic
orbit of (\ref{eq:cle}). The Poincar\'{e} map is defined by
$\pi(v_0)=v(l,v_0)$ where $v(u,v_0)$ is the solution of
(\ref{eq:cle})  with initial condition $v(0,v_0)=v_0$.

By the standard  method of differentiation of solutions of
differential equations with respect to initial conditions to
obtain the Taylor expansion of $\pi$,   it follows that
$dv/dv_0(u)=\sqrt{\frac{a_0}{a(u)}}$, where $a_0=a(0)>0$ for
definiteness.  Therefore $\pi^\prime(0)=1$.

Differentiation shows that $q(u)={\frac{d^2v}{dv_0^2}}(u,0)$
satisfies the following linear differential equation:

$$-\frac 12\frac{a_0 a^\prime b(u)}{a(u)^2} +\frac 32\frac{ a_0a^\prime
k^3}{ a(u)^2}
-\frac 16\frac{ a_0 a^\prime k_g }{a(u)}+\frac 13 \frac {a_0
b^\prime}{a(u)}   +a(u)q^\prime +\frac 12a^\prime  q(u)=0.$$

Integration of the equation above leads to

$$q(u)=\frac {a_0}6 \frac{1}{\sqrt{a(u)}} \int_0^u
[\frac{a^\prime}{a^{5/2}}(3b(u)-9k^3) +
\frac{1}{a^{3/2}}(a^\prime k_g-2b^\prime)]du$$

To obtain the second derivative of $\pi$, integration by parts
gives:

$$\pi^{\prime\prime}(0)=q(l)= \frac 16 \sqrt{a_0} \int_0^l  k_g
\frac{a^\prime }{a^{3/2}}du. $$
This ends the proof of item i).

 To proceed consider the case where
$a(u)$ has a transversal zero say at  $u=0$. Therefore the
differential equation of curvature lines, equation \ref{eq:cle},
near $(0,0)$ can be written as:

$${\mathcal F}(u,v,[du:dv]= [-v+h.o.t]dv^2+[2u+a_1 v+
h.o.t]+[v+h.o.t]du^2=0,$$
where $a_1= ( b(0)-3k^3)/a^\prime(0)$. Consider the Lie-Cartan
line field $X=({\mathcal F}_p, p{\mathcal F}_p, -({\mathcal
F}_u+p{\mathcal F}_v))$, where $p=dv/du$. The field $X$ has three
real  singular points along the projective line - axis $p$- and
its are defined by the polynomial equation $p(p^2-a_1p-3)=0$. The
linear part of $X$ along the axis $p$ is given by:

$$\aligned DX(0,0,p)= & \left( \begin{matrix}  2 & a_1-2 p  &0\\ 2 p &
( a_1-2 p) p&0\\ 0 &0&  3
p^2-2a_1p-3
\end{matrix}\right)\endaligned$$

The eigenvalues of $DX(0,0,p)$ are $\lambda_1= 2-2 p^2+a_1p$ and
$\lambda_2= 3p^2-2a_1p-3$.

At $p=0$ it follows that $\lambda_1 \lambda_2=-6<0$ and so
$(0,0,0)$ is a hyperbolic saddle point of $X$. Evaluating
$\lambda_1 $ and $\lambda_2$ at the other singular points
$(0,0,p_i)$ of $X$ it is obtained that $\lambda_1=-(1+p_i^2)$ and
$ \lambda_2= 3+p_i^2$. Therefore these two singular points of $X$
are also hyperbolic saddle points. Near the projective line the
phase portrait of $X$ is as illustrated in  Fig. \ref{fig:p}.
Blowing down the phase portrait of $X$ it follows that the
singular point is equivalent to a Darbouxian umbilic point of type
$D_3$.

\begin{figure}[htbp]
\begin{center}
\includegraphics[angle=0, width=4.5cm]{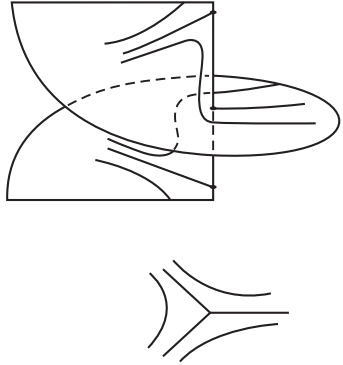}
\end{center}
\caption{ \label{fig:p} Resolution of $X$ near the projective
line}
\end{figure}
\end{proof}

\begin{remark} By properly   choosing  functions $k_g$
and $a(u)=2{\mathcal H}_v$  it is possible to construct explicit
examples of the
 spiraling behavior in item i) of Proposition \ref{prop:4}.
\end{remark}

\begin{proposition} \label{prop:3}
A closed regular curve $c:[0,l]\to \mathbb R^3$ parametrized by
arc length $u$ is a curve of umbilic points of a regular surface
containing $c$ if and only if $\int_0^l \tau(u)du\in 2k\pi, \;\;
k\in \mathbb Z$.
\end{proposition}

\begin{proof} The Frenet frame of $c$ is given by
$\{t, n, b\}$ and the following equations holds:
$$\aligned t^\prime =& k n\\
n^\prime =& -kt +\tau b\\
b^\prime =& -\tau n\endaligned $$

For any regular surface containing $c$ the  Darboux frame
$\{t,N\wedge t, N\}$ and the Frenet frame $\{t, n, b\}$  are
related by:
$$\aligned N=&\cos\theta(u) n+\sin\theta(u) b\\
N\wedge t=& -\sin\theta(u) n+ \cos\theta(u) b. \endaligned $$
Direct calculation shows that $\;k_n(u)=k(u)\cos\theta(u)$, \; $
k_g(u)=-k(u)\sin\theta(u)$\; and
$\tau_g(u)=-(\theta^\prime(u)+\tau(u))$.

Supposing that $c$ is a curve of umbilic points it follows, see
proposition \ref{prop:ca}, that $\tau_g=0$ and
$k_n(c(u),T)=k_n(c(u),N\wedge T)=k(u)\cos\theta(u)$.

Therefore in order to obtain a regular surface containing $c$
(closed curve) it is necessary that $\int_0^l\tau(u)du$ be an
integer  multiple of $2\pi$. Clearly, from  equation (\ref{eq:ca})
this condition is also sufficient.
\end{proof}

The following corollary is the first case discussed in
(Carath\'eodory 1935).

\begin{corollary} Any regular spherical closed curve $c$ is a curve of
umbilic points  of a regular surface which contains $c$.
\end{corollary}
\begin{proof} For  any closed spherical curve it follows that
$\tau_g=0$ and therefore  $\int_0^l\tau(u)du=-\int_0^l
\theta^\prime(u)du=0.$
\end{proof}

\section{Concluding Remarks} \label{sec:cr}

The interest on the structure of principal lines in a neighborhood
of a continuum of umbilic points, forming a curve,  in an analytic
surface goes back to the work of Ca\-ra\-th\'eo\-dory
(Ca\-ra\-th\'eo\-dory 1935).
  For previous results related to this subject,
Carath\'eodory refers in his paper to the books of Monge (Monge
1850)
 and Dupin  (Dupin 1813), reliquae not found by the authors.

In the present work has been  carried out   an independent,
explicit,
 self sufficient and constructive study  of umbilic curves on smooth
surfaces,
 that holds also for regular closed curves.
Here only the simplest, least degenerate cases, have been
  considered.
A partial comparison of the  results of this paper with those of
(Carath\'eodory 1935)
 is attempted  below.

In (Carath\'eodory 1935)  Carath\'eodory  established that any
local analytic regular arc of curve in  $\mathbb R^3$ is a curve
of umbilic points of a piece of  analytic surface.
 Proposition \ref{prop:3}, which also holds in the analytic case, gives
  a global independent  version.

 In some cases Carath\'eodory  also determined the behavior of principal
  lines near the curve of umbilic points. The transversal case in
Proposition
   \ref{prop:1} can be found in his paper. The tangential case is not
there.

The cases studied in Proposition \ref{prop:11} are not treated in
(Carath\'eodory 1935).

Concerning Proposition \ref{prop:4}, case $i)$, which gives an
explicit criterion for quadratic spiraling approach to the umbilic
closed curve, it seems that Carath\'eodory was aware of the
possibility of spiraling, but gave no criterion or example for
this situation. Case $ii)$ in  Proposition \ref{prop:4} was not
treated in his paper.

In
 (Carath\'eodory 1935) the focal surfaces --caustics--
received great attention.  In fact Ca\-ra\-th\'eodory starts his
analysis with the focal surfaces from which he obtains the surface
with an umbilic curve. In the present direct approach, the focal
surfaces can be obtained from their  standard expression
$$ \alpha + r_i N,$$
where $r_i = (k_i)^{-1}, i= 1, \, 2$ are the curvature radii
defined  in terms of the principal curvatures $k_{2,1} = {\mathcal
H }\pm  \sqrt{{\mathcal H}^2 -{\mathcal K}}$,  expressed in
function of the Mean and Gaussian curvatures given in  equation
\ref{eq:hk}.

Additional analysis -- complemented with some plotting -- must be
done to fully grasp the diversity of focal surfaces possible for
surfaces with an umbilic curve exhibiting the several points
studied here. In (Carath\'eodory 1935)
 only the first
cases in Propositions \ref{prop:1} and \ref{prop:4} seem to have
been considered.

 \section*{References}

 \noindent {\sc C. Carath\'eodory.} 1935.
Einfach Bemrkungen uber Nabelpunktscurven.  Lecture at Breslau,
Complete Works. Vol 5: 26-30. \vskip.2cm

 \noindent {\sc G.
Darboux.} 1896. \emph{Le\c cons sur la Th\'eorie des Surfaces,
vol. IV Sur la forme des lignes de courbure dans la voisinage d'un
ombilic , Note 07},
  Gauthier Villars,  Paris.

\vskip.2cm

\noindent {\sc Ch. Dupin}. 1813.
    D\'eveloppements de  G{\'e}ometrie, Paris,  Courcier.

\vskip.2cm

\noindent
 {\sc C. Gutierrez} and {\sc J. Sotomayor}. 1982. \emph
{Structural Stable Configurations of  Lines of Principal
Curvature},
  Asterisque {\bf 98-99}: 185-215.
\vskip.2cm

\noindent
 { \sc C. Gutierrez} and {\sc J. Sotomayor}. 1991.
{\emph  Lines of Curvature and Umbilic Points on Surfaces},
Brazilian $18^{\hbox{th}}$ Math. Coll., IMPA,
 Reprinted as \emph{ Structurally Configurations of Lines of Curvature
 and Umbilic Points on Surfaces, Monografias del IMCA}.
\vskip.2cm

\noindent
 {\sc C. Gutierrez} and {\sc J. Sotomayor}. 1998. \emph{Lines of
Curvature,  Umbilical Points and Carath\'eodory Conjecture},
 Resenhas IME-USP,  {\bf 03}: 291-322
\vskip.2cm

\noindent
  {\sc R. Garcia } and {\sc J. Sotomayor}. 1997.
 \emph{Structural stability of parabolic points and periodic asymptotic
lines},
   Matem\'atica Contempor\^anea, {\bf  12}: 83-102.

\vskip.2cm

\noindent
  {\sc R. Garcia } and {\sc J. Sotomayor}. 2000.  \emph{Lines of
Axial  Curvature on Surfaces Immersed in ${\mathbb R}^4$},
 Differential Geometry and its Applications  {\bf  12}:
253-269. \vskip.2cm

\noindent {\sc R. Garcia }, {\sc C. Gutierrez} and {\sc J.
Sotomayor}. 2000. \emph{Lines of Principal  Curvature around
Umbilics and Whitney Umbrellas},
 Tohoku Math. Journal, {\bf  52}: pp. 163-172.
\vskip.2cm

 \noindent
 {\sc R. Garcia }, {\sc C. Gutierrez} and {\sc J.
Sotomayor}. 2004. \emph{ Bifurcations of Umbilic Points and
related Principal Cycles},
 Journal of Dynamics and Differential Equations, {\bf  15}.

   \vskip.2cm

\noindent
  {\sc E. Hopf}. 1979.  Differential Geometry in the
 Large, Lect. Notes in Math. {\bf 1000},  Springer Verlag.
\vskip.2cm

\noindent {\sc G. Monge}. 1850. Application de l'Analyse  \` a  la
G{\'e}om{\'e}trie, Bachalier, Paris,   Reimpression by University
Microfilms International, Ann Arbor, MI.

\vskip.2cm

\noindent { \sc M. Spivak}. 1980.
 A Comprehensive Introduction to  Differential Geometry, Vol. III,
  Publish
or Perish, Berkeley.

\end{document}